\documentclass[11pt, a4paper]{article}

\usepackage[centertags]{amsmath} 
\usepackage{amsthm,amssymb,amscd,amsfonts}
\usepackage{mathtools}
\usepackage[all]{xy}

\usepackage[english]{babel}
\usepackage[dvipdfm,bookmarks,breaklinks]{hyperref}

\usepackage{breqn}

\theoremstyle{plain}

\newtheorem{thm}{Theorem}[section]
\newtheorem{cor}[thm]{Corollary}
\newtheorem{fact}[thm]{Fact} 
 
\newtheorem{lemma}[thm]{Lemma}
\newtheorem{prop}[thm]{Proposition} 
\newtheorem{defn}[thm]{Definition}
\theoremstyle{remark}
\newtheorem*{rem}{Remark}

\newcommand{\Aa}{\mathcal{A}}
 
\newcommand{\Ll}{\mathcal{L}} 
 
\newcommand{\Oo}{\mathcal{O}}
\newcommand{\PP}{\mathbb{P}}
\newcommand{\R}{\mathbb{R}} 

\newcommand{\Z}{\mathbb{Z}}
\newcommand{\N}{\mathbb{N}} 
\newcommand{\C}{\mathbb{C}}
\newcommand{\G}{\mathbb{G}} 

\newcommand{\Q}{\mathbb{Q}}

\newcommand{\Hh}{\mathcal{H}}

\newcommand{\m}{\mathfrak{m}}

\newcommand{\Sets}{\operatorname{Sets}}
\newcommand{\RigSp}{\operatorname{RigSp}}
\newcommand{\Sch}{\operatorname{Sch}}

\newcommand{\rk}{\operatorname{rk}}

\newcommand{\Spec}{\operatorname{Spec}}
\newcommand{\Pic}{\operatorname{Pic}}

\newcommand{\alg}{\operatorname{alg}}

\newcommand{\GL}{\operatorname{GL}}

\newcommand{\ACVF}{\operatorname{ACVF}}
\newcommand{\ACVFLR}{\ACVF_{\operatorname{LR}}}
\newcommand{\VF}{\operatorname{VF}}
\newcommand{\RV}{\operatorname{RV}}
\newcommand{\rv}{\operatorname{rv}}
\newcommand{\fin}{\operatorname{fin}}
\newcommand{\trop}{\operatorname{trop}}
\newcommand{\RES}{\operatorname{RES}}
\newcommand{\LL}{\mathbb{L}}
\newcommand{\Id}{\operatorname{Id}}
\newcommand{\Vol}{\operatorname{Vol}}
\newcommand{\Var}{\operatorname{Var}}
\newcommand{\Hilb}{\operatorname{Hilb}}

\newcommand{\eps}{\varepsilon}

\newcommand{\Hom}{\operatorname{Hom}}

\newcommand{\End}{\operatorname{End}}

\newcommand{\PGL}{\operatorname{PGL}}

\newcommand{\id}{\operatorname{id}}

\newcommand{\la}{\langle} \newcommand{\ra}{\rangle}
\newcommand{\wo}{\setminus}

\newcommand{\suchthat}[2]{\{\ #1\ \mid\ #2\ \}}

\newcommand\blfootnote[1]{%
  \begingroup
  \renewcommand\thefootnote{}\footnote{#1}%
  \addtocounter{footnote}{-1}%
  \endgroup
}

\begin{document}

\bibliographystyle{alpha}

\begin{center}
  \Large \textbf{Motivic volume of families of polarized
    rigid-analytic tori}\\[2ex]
\end{center}

\begin{center}
  {\normalsize Dmitry Sustretov 
}

\blfootnote{version of \today}
{\small 
\hspace{0.15\linewidth}
\begin{minipage}[t]{0.85\linewidth}
  \begin{center} {\bf Abstract}
  \end{center}
  Let $k$ be a non-Archimedean rational valued field.  We construct
  the moduli space of linearly rigidified polarized analytic tori over
  $k$ that admit rigid-analytic uniformization by an algebraic torus
  and observe that it is in definable rigid subanalytic bijection with
  a $\PGL_N$-bundle over a polyhedral domain in an algebraic torus. We
  use this observation to prove that the Hrushovski-Kazhdan motivic
  volume of a non-Archimedean semi-algebraic family of Abelian
  varieties admitting such a uniformization fibrewise vanishes. This
  question is motivated by the conjectural geometric interpretation of
  tropical refined multiplicities of Block and Goetsche proposed by
  Nicaise, Payne and Schroeter.
\end{minipage}

}
\end{center}

\tableofcontents

\section{Introduction}

Nicaise, Payne and Schroeter propose in their paper \cite{nps} an
approach to geometric interpretation of tropical refined Severi
degrees of Block and Goettsche \cite{block2016refined}. They
conjecture that the refined tropical multiplicity equals the
$\chi_y$-genus of the non-Archimedean semi-algebraic subset of the
universal family of compactified Jacobians over the moduli space of
stable curves of fixed genus that tropicalize to the given tropical
curve, and prove the conjecture in genus 1 for curves with a single
node. The $\chi_y$-genus is assigned to a semi-algebraic set with the
help morphism from the Grothendieck ring of (non-Archimedean)
semi-algebraic subsets to the $K$-ring of varieties over the residue
field, constructed using the theory of Hrushovski and Kazhdan
\cite{hk}. The image of a particular semi-algebraic set under this
morphism is called its motivic volume.

In view of conjectures proposed in \cite{nps} it is natural to try to
find the contribution of the semi-algebraic families of Jacobians of
smooth Mumford curves to the tropical multiplicities. More generally,
one considers semi-algebraic family of totally degenerate Abelian
varieties. In this note we prove that the motivic volume of the total
spaces of such a family is zero.

The computation of motivic volume of a family is hindered in general
due to the lack of an appropriate Fubini-type statement. In the
situation of interest, totally degenerate Abelian variety is a
quotient of an algebraic torus by a lattice, i.e. an analytic
torus. It is therefore in a definable bijection with a domain of the
form $\trop^{-1}(\Delta)$ where $\trop: \G_m^n(K) \to \R^n$ is the
coordinate-wise application of the map $-\log |\cdot|$, and $\Delta$
is an polyhedron in $\R^n$ with some of its faces removed. The motivic
volume of such semi-algebraic sets can be directly
computed. Unfortunately, in order to compute the motivic volume of a
family of such tori, the uniformization by an algebraic torus should
be uniform.

To this end we consider the moduli space of linearly rigidified
polarized analytic tori and show that the uniformization map is
locally definable in the expansion of the language of algebraically
closed valued fields with rigid subanalytic functions of Lipschitz and
Robinson \cite{lipshitz1993rigid}. We then use the invariance of
motivic volume under bijections definable in this expansion to deduce
the vanishing of the motivic volume.

The main result is Theorem~\ref{volvanish}. Section~\ref{prelim}
collects the necessary auxiliary statements about analytic tori,
polarizations, moduli of polarized Abelian varieties and
Lipshitz-Robinson rigid subanalytic functions.

\textbf{Acknowledgements}. I would like to express my gratitude to
Johannes Nicaise who have asked me the question that lead to the
appearance of this note, and to Max Planck Institute for Mathematics
for excellent working conditions.

\section{Preliminaries}
\label{prelim}

\subsection{Analytic tori, Abelian varieties and polarizations}

In order to establish notation we recall basic facts about
polarizations on Abelian schemes; we then survey the facts about
polarizations of rigid-analytic tori loosely following
\cite[Section~2.7]{lutke}, see also
\cite[Chapter~6]{fresnel2012rigid}. It is helpful to keep in mind that
the theory largely parallels the one in the complex case (see, for
example, \cite[Chapters~4 and 8]{birkenhake2013complex}).

\begin{defn}[Abelian scheme]
  A group scheme $A\to S$ is called an Abelian scheme if it is smooth, 
  proper, and its geometric fibres are connected.
\end{defn}

If $A\to S$ is an Abelian $S$-scheme, then $\Pic(A/S)$ is a smooth
proper group $S$-scheme that represents the Picard functor. Let
$\Pic^\tau(A/S)$ be its open subscheme whose geometric points
correspond to invertible sheaves that are algebraically equivalent to
zero. This scheme is smooth and projective over $S$ and its geometric
fibres are reduced and connected. The Abelian scheme $\Pic^\tau(A/S)$
is called the scheme \emph{dual to $A$} and is denoted $\hat{A}$. A 

The universal line bundle on $A \times \hat{A}$ is called
\emph{Poincar\'e line bundle} and is denoted $P_{A \times \hat{A}}$. 

Let $L$ be an arbitrary line bundle on $A$, and let $\mu: A \times_S A
\to A$ be the multiplication map. Then the line bundle
$$
\mu^* L \otimes (p_1^*L)^{-1} \otimes (p_2^*L)^{-1}
$$
can be considered as a line bundle over $X$ via the projection
$p_1: A \times_S A \to A$, and so by the definition of $\Pic$ gives
rise to the mapping $\omega_L: A \to \Pic(A/S)$. If $e: S \to A$ is
an identity then $\omega_L \circ e$ is the identity of
$\Pic(A/S)$. Since the fibres of $\Pic(A/S)$ are connected, the
morphism $\omega_L$ factors through $\Pic^\tau(A/S)$.

Recall that a morphism of Abelian varieties over a field is called an
isogeny if it is surjective with finite kernel; a morphism of Abelian
schemes is a morphism that induces isogenies on geometric fibres.

\begin{fact}
  The construction above induces a homomorphism
  $\Pic(A/S) \to \Hom_{\Z}(A,\hat{A})$. The map $\omega_L$
  is an isogeny if and only if $L$ is ample.
\end{fact}

\begin{defn}[Polarisation]
  A polarization of an Abelian $S$-variety $A$ is a morphism $\varphi:
  A \to \hat{A}$  such that for all geometric fibres the induced
  morphism $\varphi_s: A_s \to \hat{A}_s$ is an isogeny of the form
  $\omega_{L_s}$ for some ample line bundle $L$ on $A$.
\end{defn}

Fix a non-Archimedean valued field $k$, let $M$ be a free Abelian
group of rank $n$, denote by $M'$ its dual $\Hom(M,\Z)$, and let
$T := \Spec k[M']$. Denote by $\trop: T \to \R^n$ the
coordinatewise valuation map:
$$
\trop(x_1, \ldots, x_n) = (-\log |x_1|, \ldots, -\log |x_n|)
$$
A torsion-free subgroup $\Lambda \subset T$ is called a \emph{lattice}
if $\trop$ induces an isomorphism between $\Lambda$ and a discrete
subgroup $\trop(\Lambda)$ of the additive group $\R^n$.

Lattices $M \to T$ are in natural bijective correspondences with the
lattices $M' \hookrightarrow T'$. Indeed, regarding $M'$ as $\Hom(T,
\G_m)$ and $T'$ as $\Hom(M,\G_m)$ the embedding $M' \hookrightarrow T$
is given by the restriction map $\Hom(T,\G_m) \to \Hom(M,\G_m)$. 

\begin{fact}[Proposition~2.7.5 \cite{lutke}, Non-Archimedean Appel-Humbert theorem]
  The set of isomorphism classes of line bundles on $T/M$ is in
  bijective correspondence with pairs $(\lambda, r)$, where
  $\lambda: M \to M'$ is a homomorphism and $r: M \to \G_m$ subject to
  the condition
  $$
  \lambda(m_1)(m_2) = \dfrac{r(m_1 + m_2)}{r(m_1)r(m_2)}
  $$
  $\lambda$ is trivial if and only if $L_{(\lambda,r)} \in \Pic^0$,
  moreover, $\Pic^0$ consists of groups of translation-invariant line
  bundles. 
\end{fact}

One observes that the function
$Z: M \to H^0(T,\Oo_T^\times) = r(m)\lambda(m)$, called the
\emph{automorphy factor} is a group cohomology 1-cocycle,
$Z \in H^1(M, H^0(T,\Oo_T^\times))$. A function
$f \in H^0(T,\Oo_T^\times)$ is called \emph{theta function with
  respect to the automorphy factor $Z$} if 
$$
f(m\cdot x) = Z_m f(x)
$$

\begin{defn}[Polarization of an analytic torus]
  A polarization of the analytic torus $T/M$ is an injective map
  $\lambda: M \to M'$ such that the bilinear map
  $$
  \la \cdot, \cdot \ra: M \times M \to K^\times\ \ \ \ \la m_1, m_2 \ra =
  \lambda(m_1)(m_2)
  $$
  is symmetric and positive definite, that is, for any $m\in M$, $-
  \log |\la
  m,m \ra| > 0$.  
\end{defn}

\begin{rem}
  Let $M \hookrightarrow T$ be a lattice. If $\lambda: M \to M'$ is a
  homomorphism of groups then it induces a morphism of tori
  $\varphi_\lambda: T \to T'$. If $\lambda$ induces a symmetric and
  non-degenerate map $\la \cdot, \cdot \ra$ then
  $\varphi_\lambda(M) \subset M'$ and so the morphism
  $\varphi_\lambda: T/M \to T'/M'$ is well-defined. If $\lambda$
  defines a polarization then $\varphi_\lambda = \varphi_L$ for an
  ample line bundle $L$.
\end{rem}

\begin{fact}[Theorem~2.7.12, \cite{lutke}]
  A line bundle is ample if and only if $\lambda$ defines a
  polarization. The global sections of a line bundle $L$ are given by
  theta functions with respect to the automorphy factor $Z$.
\end{fact}

\begin{fact}[Lemma~6.5.4, \cite{fresnel2012rigid}]
  \label{triple-emb}
  If $L$ is an ample line bundle on $T/M$ and
  $\theta_0, \ldots, \theta_n$ is a basis of $H^0(T/M,L^3)$ then $x
  \mapsto (\theta_0(x) : \ldots : \theta_n(x))$ defines a closed
  embedding $T/M \hookrightarrow \PP(H^0(T/M,L^3))$.
\end{fact}

\begin{fact}[Propositions~6.10 and 6.13, \cite{mumford-git}]
  \label{polar-double}
  Let $\omega: A \to \hat{A}$ be a polarization of an Abelian variety,
  let $L = (\id \times \omega)^* P_{A \times }$ and let $\omega'$ be
  the polarization induced by $L$. Then
  \begin{itemize}
  \item $\omega' = 2\omega$
  \item $(\dim H^0(A, L))^2 = \deg \omega$
  \end{itemize}
\end{fact}

Denote $\rk \lambda=\#(M'/\lambda(M))$. The following fact easily
follows from the automorphy equation.
\begin{fact}
  \label{sections}
  Theta functions have the form
  $f(x) = \sum_{\chi \in M'} a_\chi \chi$, and are determined
  by coefficients $a_{\chi_1}, \ldots, a_{\chi_n}$ where
  $\chi_1, \ldots, \chi_n$ are representatives of $M'/\lambda(M)$. In
  particular, $\dim H^0(T/M, L) = \rk \lambda$.
\end{fact}

If $L$ is an arbitrary line bundle, define
$$
\varphi_L: T/M \to T'/M'\ \ \ \ \varphi_L(a) = t^*_a L \otimes L^{-1}
$$
where $t_a: T/M \to T/M, t_a(x) = x+a$ for any $a\in T/M$. Clearly,
the line bundle $t^*_a L \otimes L^{-1}$ is translation-invariant, so
the map is well-defined. One can show that $\varphi_L$ is an analytic
homomorphism of groups.

\begin{fact}
  \label{polar-deg}
  The line bundle $L$ is ample if and only if $\varphi_L$ is
  surjective. The degree of $\varphi_L$ is of size $d^2$ where $d$ is
  the degree of $L$.
\end{fact}

\begin{fact}[Riemann-Roch on an Abelian variety]
  \label{rr-abelian}
  Let $L$ be a positive line bundle on an Abelian variety of dimension
  $g$, then
  \begin{align*}
    \chi(L) & = L^g/g!\\
    \chi(L)^2 & = \deg \varphi_L\\
  \end{align*}
  Consequently, the Hilbert polynomial of an Abelian variety endowed with
  polarization $\varphi$ of degree $d$ with respect to
  $L_\varphi^{\otimes 3}$ is $P(x) = x^g d$.
\end{fact}

\subsection{Motivic volume}

Hrushovski-Kazhdan motivic integration theory
\cite{hk} provides a way to express
non-Archimedean semi-algebraic subsets of algebraic varieties over a
valued field $K$ as unions of semi-algebraic sets of two particular
kinds. The first one is related to the geometry of integral polhedra,
and the second one is related to algebraic varieties over the residue
field. 
 
Formally, the theory is formulated in the context of model theory of
algebraically closed valued fields. Let $K$ be such a field, then one
considers several \emph{sorts}: the valued field sort $\VF$, the
residue-value sort $\RV$ and the value group sort $\Gamma$. 

Let $\Oo \subset K$ be the value ring with the maximal ideal $\m$.
Consider the exact sequence of groups
$$
1 \to \Oo^\times/(1+\m) \to K^\times/(1+\m) \to \Gamma \to 0
$$
The middle term is called $\RV$ and is made into a sort with the
structure of the multiplicative group, and two inter-sort projection
maps: $\rv: \VF \wo \{0\} \to \RV$, and $v_{\rv}: \RV \to \Gamma$. 

After fixing some base field $K_0$, one associates the following
categories of definable sets to the sorts $\VF, \RV$ and $\Gamma$.

\begin{defn}
  The category $\VF[n]$ is defined to be the category of definable
  subsets of of $n$-dimensional varietis over $K_0$.

  The category $\RV[n]$ is defined to be the category of pairs $(X,f)$
  where $X$ is a definable set and $f: X \to \RV^n$ is a definable map
  with finite fibres.

  The category $\Gamma[n]$ is defined to be the category of Boolean
  combinations of subsets of $\Gamma^n$ defined by inequalities and
  equalities with integral coefficients and with parameters in
  $\Gamma(K_0)$. 
  
  $\Gamma^{\fin}[n]$ is the full subcategory $\Gamma[n]$ of definable
  finite subsets.

  $\RES[n]$ is the full subcategory of $\RV[n]$ which consists of
  definable sets which project to finite definable subsets of
  $\Gamma^n$ via $\trop$.
\end{defn}

\begin{defn}
  If $A$ is a category of definable sets then we denote by $K_+(A)$
  the semi-ring generated by definable subsets in $A$ module the
  relations 
  \begin{itemize}
  \item $[A] = [B]$ if there exists a definable bijection between $A$
    and $B$,
  \item $[C] = [A] + [B]$ if $C = A \sqcup B$.
  \end{itemize}
\end{defn}

In case $K_0=k((t))$ there exists a canonical isomorphism between a
quotient $!K(\RES)$ of the ring $K(\RES)$ and the equivariant
Grothendieck ring $K^{\hat \mu}_0(\Var_k)$. Let
$\theta: K_0(\RES) \to K_0\Var$ be the composition of the quotient map
with this canonical isomorphism and the forgetful morphism
$K^{\hat \mu}_0(\Var_{k_0}) \to K_0(\Var_k)$ (see Seciton~4 of
\cite{hrulo}).

Denote
$K_+\VF=\oplus_n K_+(\VF[n]), K_+\RV[\leq n]=\oplus_{l \leq n}
K_+\RV[l]$, and define the morphisms
$$
\begin{array}{cll}
  & \LL: K_+(\RV[n]) \to K_+(\VF[n]), &\ [(X,f)] \mapsto [\VF^n
                                        \times_{\rv,F} X]\\
  & \LL: K_+(\Gamma[n]) \to K_+(\VF[n]),&\   [\Delta] \mapsto [\trop^{-1}(\Delta)] \\
\end{array}
$$

The motivic integration theory of Hrushovski and Kazhdan
\cite{hk} (in the non-measured case) rests on
two main statements: that the natural morphism
$$
\oplus_{l+m=n} K_+(\Gamma[l]) \otimes_{K_+\Gamma^{\fin}} K_+(\RES[m]) \to
K_+(\RV[n])
$$
is an isomorphism and that the morphism
$$
\LL: \oplus_n K_+(\RV[\leq n]) \to K_+(\VF)
$$
is surjective. The kernel $I_{sp}$ of the latter can be explicitly
described. The theory is developed in an axiomatic setting that
depends only on the category $\RV$ (the corresponding notion is called
$V$-minimality).

Consider modified Euler characteristic 
$$
\chi': \oplus_n K_+(\Gamma[n]) \to \Z, \chi'([\Delta]) = \lim_{n\to \infty}
\chi(\Delta \cap [-l,l]^n)
$$
where $\chi$ is the usual o-minimal Euler characteristic (which yet
again coincides with the usual Euler characteristic when $\Gamma \cong
\R$).

Define the morphism $\Vol: K_0(\VF[n]) \to K_0(\Var_{k_0})$
$$
\Vol (\LL^{-1}([X] \otimes [\Delta]) = \theta([X])\cdot \chi'(\Delta)(\LL-1)^n
$$
is well-defined because $\Id \otimes \chi'$ is trivial on
$I_{sp}$. The destination of the morphism can be identified with
$K_+(\Var_k)$ if $K$ is algebraically closed.

\subsection{Lipschitz-Robinson rigid subanalytic functions}

Let $K$ be an algebraically closed, complete, non-Archimedean normed
field. Let $R = \suchthat{ x \in K}{ |x| \leq 1}$, let
$\m=\suchthat{ x \in K}{ |x| < 1}$, and let $k=R/\m$. Define the norm on the ring
$K[[x, \rho]]$ as follows:
$$
|\sum a_{ij} x^i \rho^j| = \sup |a_{ij}|
$$

Let $R_0 \subset R$ be a maximal discrete valuation ring contained in
$R$ with prime $\pi \in \m$ such that $0 < |\pi| < 1$ and
$R_0/(\pi) \cong k$. For any sequence $(a_i)$ with $a_i \in R$ and
such that $|a_i| \to 0$ let $\widehat{R_0[\{a_i, i \in \N\}]}$ be the
completion of $R_0[\{a_i, i \in \N\}]$ with respect to the norm on $K$
and define
$$
R_0\{a_i\}\la x \ra = \widehat{R_0[\{a_i, i \in \N\}]}\la x \ra
$$
Let $R_0\{a_i\}\la x \ra[[\rho]]$ be the ring of formal power series
over $R_0\{a_i\}\la x \ra$. Define 
$$
S\{a_i\}\la x \ra[[\rho]] = \suchthat{\pi^{-\alpha} f}{f \in
  R_0\{a_i\}\la x \ra[[\rho]]}
$$
and 
$$
K\la x \ra[[\rho]]_s = \cup_{\{a_i\}} S\{a_i\}\la x \ra [[\rho]]
\subset K[[x,\rho]]
$$
The elements of this ring define analytic functions $R \times \m \to
K$ which are well-behaved. For example, these functions have finitely
many zeroes on $R \times \m$.

A \emph{rigid subanalytic function} is a function definable in the
expansion of $K$ with graphs of functions from $K\la x \ra[[\rho]]_s$ (Lipschitz and Robinson \cite{lipshitz1993rigid}).

\begin{prop}
  \label{def-functions}
  Let $S \subset X$ be a semi-algebraic subset of an algebraic variety
  $X$, and assume that $S$ is a finite union of rational and
  semi-rational domains.  An analytic function on a proper
  semi-algebraic subset of an algebraic variety is definable in the
  language $\ACVFLR$.
\end{prop}

\begin{proof}
  Follows from the fact that analytifications of affine varieties can
  be covered by affinoid domains, that functions analytic on rational
  and semi-rational subdomains of affinoid domains are definable, and
  that semi-algebraic domains are contained in finite unions of
  rational and semi-rational domains.
\end{proof}

\begin{cor}
  \label{fund-domain}
  Let $T=\G_m^g$ be a torus, and assume that discrete group $G$ acts
  on $T$ so that fundamental domain $U \subset T$ is
  semi-algebraic. Let $f$ be a meromorphic funtion on $T/G$, and let
  $p: T \to T/G$ be the quotient map. Then restriction of $p \circ f$
  to the fundamental domain is definable in $\ACVFLR$. 
\end{cor}

As was remarked in the previous section, the motivic integration
theory of \cite{hk} can be carried out verbatim in any expansion of
the theory of algebraically closed valued fields as long as the
structure induced on the sort $\RV$ is unchanged. In particular the
following is true.

\begin{fact}[Lemma~3.33, \cite{hk}]
  \label{acvflr}
  Let $X,X'$ be semi-algebraic subsets alrgebraic varieties over an
  algebraically closed valued field. If there exists a bijection
  $X \xrightarrow{\sim} X'$ definable in $\ACVFLR$ then $[X] = [X']$
  in $K_0(\VF)$.
\end{fact}

\section{Motivic volume of a family of polarized analytic tori}

In this section $k$ is a complete rational valued field, i.e. a field
complete with respect to a non-Archimedean norm and such that the
image of the map $\log |\cdot|: k^\times \to \R$ is contained in $\Q$,
for example, $k$ can be a discretely valued field, or its algebraic
closure, such as the field of Laurent series $\C((t))$ or the field of
Puiseux series $\C((t))^{\alg}$. Denote the residue field $\bar k$.

\subsection{The moduli space of linearly rigidified polarized analytic
  tori}

\begin{defn}[Uniformized analytic tori]
  By a family of uniformized analytic tori we will understand
  \begin{itemize}
  \item flat morphism of rigid analytic spaces $\pi: A \to S$
  \item an action of $\Z^g$ on $\G_m^g \times S$ by shifts so that
    $\Z^g \hookrightarrow (\G_m^g)_s$ is a lattice for all $s \in S$,
    and an $S$-isomorphism $\G_m^g \times S \cong A$.
  \end{itemize}
  Two families $A_1\to S$ and $A_2\to S$ are isomorphic if there
  exists an $S$-isomorphism $A_1 \xrightarrow{\sim} A_2$ that can be
  lifted to a $\Z^g$-equivariant isomorphism of respective covers by
  $\G_m^g \times S$.
\end{defn}

\begin{defn}[Linear rigidification]
  Let $S$ be a scheme or an analytic space, let $p: A \to S$ be an
  analytic torus or an Abelian scheme over $S$ and let
  $\varphi: A \to \hat{A}$ be a polarization of degree $d$. Then for
  any $s\in S$
  $$
  \dim (p_* \Ll_\varphi^3)_s = m := 6^g\cdot d
  $$
  An isomorphism $\PP(p_* \Ll_\varphi^3) \cong \PP(\Oo_S^m)$ is called
  a \emph{linear rigidification of $(A,\varphi)$}.
\end{defn}

Let $M, M'$ be free rank $g$ Abelian groups, and let
$T = \Spec k[M'], T'=\Spec k[M]$. Let $\lambda: M \to M'$ be an
injective homomorphism and $S$ be a rigid analytic space. Define
$\Aa^u_{g,\lambda}(S)$ to be the set of isomorphism classes of
uniformized polarized analytic tori with polarization of type
$\lambda$ and define $\Hh_{g,\lambda}^u(S)$ to be the set of
isomorphism classes of uniformized polarized analytic tori together
with a linear rigidification. This defines two functors
$$
\Aa_{g,\lambda}^u: \RigSp_k \to \Sets\qquad \Hh_{g,\lambda}^u: \RigSp_k \to \Sets
$$

Pick a distinguished basis $\eps_1, \ldots, \eps_g$ in $M$, then the
space $B_g$ of all embeddings $M \hookrightarrow T$ with the
distinguished basis be identified with the space of matrices
$E=(e_{ij})$ where $e_{ij} \in K^\times$ is the $j$-th coordinate of
the image of $i$-th basis vector.  Define the space of lattices with a distinguished basis. 
$$
\tilde{B}_g =  \suchthat{ \iota(\eps_1), \ldots,
  \iota(\eps_g) }{ \iota: M \hookrightarrow T}
$$
The group $\GL(M)$ acts on the
on  $B_g$: if $\Omega=(\omega_{ij}) \in \GL(M)$ then
$$
\Omega \cdot E = (\prod_i e_{ij}^{\omega_{ji}})
$$
and the quotient $\tilde{B_g}$ is the space of embeddings
$M \hookrightarrow T$.

We call a domain $A \subset \G_m^n$ \emph{polyhedral} if $A =
\trop^{-1}(\Delta)$ for some integral polytope $\Delta$. 

\begin{prop} 
  \label{glm-fund}
  The fundamental domain of the monomial free action of $\GL(M)$ on
  $\tilde{B}_g$ is polyhedral.
\end{prop}

\begin{proof}
  This is easily deduced from the fact that $\GL_n(\Z)$ is generated
  by diagonal matrices which have $+1$ and $-1$ entries and matrices
  of the form $\Id + E_{ij}$ where $E_{ij}$ is the elementary matrix
  that has 1 as the $ij$ entry and otherwise 0.
\end{proof}

Fix an isomorphism $i: M \cong M'$, then for any embedding $M
\hookrightarrow T$ given by the matrix $E$ the corresponding embedding
$M' \hookrightarrow T'$ is represented by $E^*$.

Define the universal tori 
$$
\tilde{Z}_g = (T \times \tilde{B}_g)/M \qquad \tilde{Z}'_g = (T'
\times \tilde{B}_g)/M'
$$
over $\tilde{B}_g$. Since the map $\varphi_\lambda: T \to T'$ is
$M$-equivariant, it descends to the quotients. Furthermore, the
action of $GL(M)$ naturally lifts from $B_g$ to $T \times B_g$ and
sends $T_s$ to $T_{\Omega s}$ is such a way that $\Omega(M_s) =
M_{\Omega s}$.

Any homomorphism $\lambda: M \to M'$ is of the form $\Lambda \circ i$;
if $\lambda$ is a polarization then $\Lambda \in \End(M)$ is
injective.

The morphism $\lambda$ induces a surjective morphism of algebraic tori
$\varphi_\lambda: T \to T'$ and for any embedding
$M\hookrightarrow T$, $\varphi_\lambda(M) \subset M' \subset T'$ and
$\varphi_\lambda|M=\lambda$. It therefore descends to the quotients:
$\varphi_\lambda: Z_g \to Z'_g$.

For any matrix $E \in B_g$, $E=(e_{ij})$ denote by $\bar{E}$ the
matrix $(-\log |e_{ij}|)$.  Define
$$
\tilde{A}_{g,\lambda} = \suchthat{ E \in B_g}{ (\Lambda E) = (\Lambda E)^*,
  \bar E > 0}
$$
Put $\tilde{Z}^u_{g,\lambda} = Z_g \times_{\tilde{B_g}}
\tilde{A}^u_{g,\lambda}$. 

By construction, each $x \in A^u_{g,\lambda}$ defines a lattice and
the fibre $(Z^u_{g,\lambda})_x$ carries the structure of a uniformized
analytic torus, and the restriction of $\varphi$ to it is a
polarization.

\begin{prop}
  \label{polarization}
  If $\iota: M \hookrightarrow T$ is represented by a matrix
  $E \in \tilde{A}_{g,\lambda}$ then the map
  $\varphi_\lambda: T/M \to T'/M'$ is a polarization.
\end{prop}

\begin{proof}
  We need to check that the form
  $\la -, - \ra: M \times M \to K^\times, \la m_1, m_2 \ra =
  \lambda(a)(b)$ is symmetric and positive definite. Indeed,
  \begin{dmath*}
    \left\la \sum_{i=1}^g a_i \eps_{i}, \sum_{j=1}^g b_j \eps_{j}
    \right\ra = \lambda(\sum_i a_i \eps_i)(\prod_j
    \iota(\eps_j)^{b_j}) = \prod_k (\prod_j e_{kj}^{b_j})^{\sum_i
      \lambda_{ik} a_i } = \prod_i \prod_j (\prod_k
    e_{kj}^{\lambda_{ik}})^{a_i b_j}
  \end{dmath*}
  which is clearly symmetric, given
  $\prod_k e_{kj}^{\lambda_{ik}} = \prod_k
  e_{ki}^{\lambda_{jk}}$. Further,
  $$
  - \log |\la \sum a_i \eps_{i}, \sum_{j=1}^g b_j \eps_{j} \ra| =
  \sum_i \sum_j  a_i b_j (\sum_k \lambda_{ik} \bar{e}_{kj}) = (\Lambda
  \bar{E} a, b)
  $$
  where $(-,-)$ is the Euclidean scalar product on $\R^n$. Therefore,
  since the matrix $\Lambda \bar{E}$ is strictly positive definite,
  the bilinear symmetric form $-\log |\la a, b \ra|$ is also positive
  definite.
\end{proof}

\begin{prop}
  \label{uni-poly}
  The space $\tilde{A}^u_{g,\lambda}$ is isomorphic to a union of
  polyhedral domains. For any polyherdal domain
  $S \subset \tilde{A}^u_{g,\lambda}$, there is a bijective morphism
  from a polyhedral domain onto
  $\tilde{Z}_{g,\lambda} \times_{\tilde{A}^u_{g,\lambda}} S$.
\end{prop}

\begin{proof}
  The first statement clearly holds for $\tilde{A}^u_{g,\id}$: the
  symmetry condition is intersection of some diagonal varieties, and
  positivity condition means that the coefficients of $\bar{e}_{ij}$
  belong to some open subset of $R^{g^2}$, which is a union of
  integral polyherdra. Notice that
  $$
  \tilde{A}^u_{g,\lambda} = \tilde{A}_{g,\id} \times_{B_g,\lambda} B_g
  $$
  For any polyhedral domain $\trop^{-1}(\Delta) \subset A^u_{u,\id}$
  the set $\trop^{-1}(\Delta) \times_{B_g,\lambda} B_g$ is polyhedral
  since $\lambda$ is a monomial morphism.

  The second statement follows from the fact that the fundamental
  domain of the action of $M$ on each fibre
  $(T \times \tilde{A}^u_{g,\lambda})_s$, is a polyhedral domain, and
  that it only depends on $\trop(s)$.
\end{proof} 

Let $T/M$ be an analytic torus with a polarization
$\varphi_\lambda: T/M \to T'/M'$ and let
$L = (\id \times \varphi_\lambda)^* P$ where $P$ is the Poincar\'e
bundle on $T/M \times T'/M'$. Then linear rigidifications of $L^3$ are
a $\PGL_N$ torsor where $N=6^g d$ and $d=\rk \lambda$ (by
Facts~\ref{polar-double} and \ref{triple-emb}). Indeed, by
Fact~\ref{sections} the sections of $L_\lambda$ are determined by
coefficients $a_{w_1}, \ldots, a_{w_N}$ of theta functions
$\sum_{\chi \in M'} a_\chi \chi$, where $w_1, \ldots, w_N$ are some
representatives of $M'/6\cdot\lambda(M)$. A linear rigidification is
uniquely determined by a choice of basis in the space of these
coefficients, up to scalar multiplication.

Acting by automorphism of the torus on the argument sends characters
of $T$ to characters. Let
$\tilde{H}^u_{g,\lambda} = \tilde{A}^u_{g,\lambda} \times \PGL_N$ and
extend the action of $\GL_g(\Z)$ from $\tilde{A}^u_{g,\lambda}$ to
$\tilde{H}^u_{g,\lambda}$ via the action of $\GL_g(\Z)$ on the basis
of theta functions by substitution:
$$
\Omega \cdot \sum_{\chi \in M'} a_\chi \chi(x) = \sum_{\chi \in M'}
a_{\Omega \chi} \chi(x)
$$

Define
$$
Z^u_{g,\Lambda} = \tilde{Z}_{g,\lambda}/\GL_g(\Z) \qquad
H^u_{g,\Lambda} =  \tilde{H}_{g,\lambda}/\GL_g(\Z) \qquad
A^u_{g,\Lambda} =  \tilde{A}_{g,\lambda}/\GL_g(\Z)
$$
The obvious map that forgets linearization makes $H^u_{g,\lambda}$
into a $\PGL_N$-bundle over $A^u_{g,\lambda}$.

\begin{lemma}
  \label{rational-stuff}
  Let $A \in M_g(\R)$ be a real matrix, and assume that there is a
  neighbourhood $U$ of $A$ in $M_g(\R)$ such that all matrices $A' \in
  U$ define bilinear forms that are strictly positive definite on
  $\Z^n \subset \R^n$. Then $A$ is positive definite. 
\end{lemma}

\begin{proof}
  First note that $x \mapsto (Ax,x)$ is positive on $\Z^n$ if and only
  if it is positive on $\Q^n$.

  Suppose $A$ is not positive definite. It cannot have negative
  eigenvalues, so assume it has an eigenvector $x$ with eigenvalue
  0. By assumption this vector has irrational coordinates. As $A'$
  varies in $U$, the eigenspace $\R\cdot x$ varies too. Clearly there
  exists an $A'$ with arbitrarily close eigenspace $V$ with eigenvalue
  zero with $V \cap \Q^n \neq \{0\}$. For such $A'$ the assumption is
  not true, and we have arrived at a contradiction.
\end{proof}

\begin{prop}
  \label{moduli-h}
  For any polarization type $\lambda: M \to M'$ and any rationally valued
  field $k$, the set of analytic tori $(Z_{g,\lambda})_x$ as $x$
    ranges in $H^u_{g}(k)$ coincides with the set
    $\Hh^u_{g,\lambda}(k)$.  
\end{prop}

\begin{proof}
  Follows from construction of $H^u_{g,\lambda}$,$Z_{g,\lambda}$
  Propositon~\ref{polarization} and Lemma~\ref{rational-stuff}.
\end{proof}

Recall that the functor $\Hh_{g,d,n}: \Sch/S \to \Sets$, defined in
Section~6 of \cite{mumford-git}, associates to a scheme $S$ the set of
isomorphism classes of linearly rigidified degree $d$ polarized
Abelian schemes over $S$ with level $n$ structures. For our purposes
we do not need to deal with the level structure and we will only
consider the functor $\Hh_{g,d,1}$ which we will denote $\Hh_{g,d}$.
Let $H_{g,d}$ be the $k$-scheme that represents the functor
$\Hh_{g,d}$.

\begin{prop}
  For any polarisation $\lambda$ of degree $d$ there exists a
  rigid-analytic embedding
  $\Hh^u_{g,\Lambda} \hookrightarrow (\Hh_{g,d})^{an}$ and a
  rigid-analytic embedding $Z_{g,\Lambda} \to Z_{g,d}^{an}$ compatible
  with projection to $H_{g,d}$
\end{prop}

\begin{proof}
  By \cite[Theorem~4.1.3]{conrad2006relative} there exists an analytic
  embedding of $\Hh^u_{g}$ into $((\Hilb_{\PP^N}/k)^{P(x)})^{an}$,
  where $P(X) = 6^g \cdot d \cdot x^g$, and of $Z^u_{g,\lambda}$ into
    $(Z_{g,d})^{an}$. For any $s \in \Hh^u_{g}$ the fibre
    $(Z_{g,\lambda})_s$ is a polarized Abelian variety and hence, by
    Proposition~7.3 of \cite{mumford-git},
    $s \in H^{an}_{g,d} \subset ((\Hilb_{\PP^N}/k)^{P(x)})^{an}$.
\end{proof}

\begin{cor}
  \label{moduli}
  Let $k$ be rationally valued.  For any $k$-variety $S$, any
  polarized Abelian scheme $A \to S$ and any semi-algebraic subset
  $U \subset S$ such that $A_s$ is mulitplicatively uniformized for
  all $s \in U$ there exists a map $U \to H^u_{g}$ such that
  $Z \times_{H^u_{g}} U \cong A \times_S U$.
\end{cor}

\subsection{Integration}

We are going to use the tropical motivic Fubini theorem of Nicaise and
Payne which we now recall.

\begin{thm}[\cite{nicaise2017tropical}]
  \label{fubini}
  Let $A \subset Y \times \G_m^n$ be a semi-algebraic subset, and let
  $\pi: Y \times \G_m^n \to \G_m^n$ be the projection map. Then there
  definable subsets $\Delta_1, \ldots, \Delta_m \subset \R^n$ and
  classes $X_1, \ldots, X_m \in K_0(\Var_{\bar k})$ such that for any
  integer $i, 1 \leq i \leq n$ and for any $\xi \in \Delta_i$,
  $\Vol((\trop\circ \pi)^{-1}(\xi)) = X_i \in K_0(\Var_{\bar k})$ and
  $$
  \Vol(A) = \sum_{i=1}^m \chi'(\Delta_i)(\LL-1)^n \cdot X_i
  $$
\end{thm}

We finally put together all the ingredients prepared so far.

\begin{thm}
  \label{volvanish}
  Let $T$ be a $k$-variety, let $\pi: A \to S$ be a Abelian scheme of
  relative dimension $g$ over $S$ and let $U \subset S$ be a
  semi-algebraic set such that $A_s$ can be uniformized by a torus for
  any $s \in U$. Then $\Vol(A \times_S U)=0$.
\end{thm}

\begin{proof}
  By Fact~\ref{acvflr} we may use maps definable in $\ACVFLR$.

  We may assume that $S$ and $T$ are connected. Pick some polarization
  on $A$, then by Corollary~\ref{moduli} there exists a map
  $T \to H_{g,d}$ for some $d$ and such that the image of $U$
  lies in $H^u_{g,\lambda}$ for some $\lambda$, $\rk \lambda = d$.

  Using Corollary~\ref{fund-domain} and Proposition~\ref{uni-poly} we
  will identify $Z^u_{g,\lambda}$ and $A^u_{g,\lambda}$ with unions of
  polyhedral domains.  It follows from Proposition~\ref{glm-fund} that
  there exists a decompostion $U = \sqcup U_i$ with $U_i$
  semi-algebraic such that restirictions of $A \times_S U_i$ to the
  fibres of the projection $H^u_{g,\lambda} \to A^u_{g,\lambda}$ are
  trivial families of tori. Consequently, $A \times_S U$ is in a
  definable bijection with a semi-algebraic set
  $Z^u_{g,\lambda} \times_{A^u_{g,\lambda},\psi} U$ for some definable
  map $\psi: U \to A^u_{g,\lambda}$.

  Let $\Sigma = \trop(\psi(U)) \subset \trop(A^u_{g,\lambda})$. Then
  $Z^u \times_{A^u_{g,\lambda}} \trop^{-1}(\Sigma) =
  \trop^{-1}(\Delta)$ for some definable subset $\Delta \subset \R^n$
  for some $n$. Let $\psi'$ be the definable bijection
  $A \to Z^u_{g,\lambda} \times U$ induced by $\psi$. Then
  $A \times_T U$ is in definable bijection with the graph
  $\Gamma_{\psi'} \subset (A \times_T U) \times \trop^{-1}(\Delta)$ of
  the map $\psi'$.

  Denote $\pi: \Delta \to \Sigma$ the natural projection, and denote
  $$
  P=\suchthat{x \in k^\times}{ |x| = 1} 
  $$
  the unit annulus. One observes that
  $$
  (\trop \circ \psi')^{-1}(\xi) = (\trop \circ \psi)^{-1}(\xi) \times
  P^g
  $$
  for $\xi \in \Delta$.

  Finally, by Theorem~\ref{fubini} there exists a decomposition
  $\Sigma = \sqcup_{i=1}^n \Sigma_i$ into definable subsets such that
  $\Vol((\trop \circ \psi)^{-1}(\xi))$ is constant for all
  $\xi \in \Sigma_i$, for each $i$, and so
  \begin{dmath*}
    \Vol(\Gamma_{\psi'}) = \sum_{i=1,\xi \in \Sigma_i}^n \Vol(P^g
    \times (\trop \circ \psi')^{-1}(\xi)) \chi'(\pi^{-1}(\Sigma_i))(\LL
    - 1)^{\dim \pi^{-1}(\Sigma_i)} =
    \sum_{i=1,\xi \in \Sigma_i}^n \Vol((\trop \circ \psi')^{-1}(\xi)))
    \chi'(\pi^{-1}(\Sigma_i)) (\LL
    - 1)^{\dim \pi^{-1}(\Sigma_i)}
  \end{dmath*}
  Here, $\chi'(\pi^{-1}(\Sigma_i)) = 0$ since
  $\chi'$ is multiplicative and fibres of $\pi$ are fundamental
  domains of a lattice, so $\chi'$ vanishes on them.
\end{proof}

\begin{cor}
  Let $C \to T$ be a family of smooth projective curves. Let
  $S \subset T$ be a semi-algebraic subset of $T$ such that $C_s$ is a
  Mumford curve for all $s \in S$, and let $J(C/T) \to T$ be the
  relative Jacobian. Then $\Vol(J(C/T) \times_T S) = 0$.
\end{cor}

\begin{proof}
  The family $J(C/T) \to T$ is a projective Abelian scheme, and its
  restriction to $S$ can be uniformized by a torus fibrewise,
  therefore, Theorem~\ref{volvanish} applies.
\end{proof}

\bibliography{tori}

\vspace{3ex}

\noindent {\sc Max Planck Institut f\"ur Mathematik\\
 Vivatsgasse 7\\
 53111 Bonn\\
 Germany\\}
{\tt sustretov@mpim-bonn.mpg.de\\}

\end{document}